\documentclass{article}

\usepackage{amssymb}
\usepackage{dsfont}

\newcommand*\dif{\mathop{}\!\mathrm{d}}
\DeclareMathSymbol{\sm}{\mathbin}{AMSa}{"39}
\usepackage{graphicx} 
\DeclareMathSymbol{\smin}{\mathbin}{AMSa}{"39}

\usepackage{mathtools}

\usepackage{float}
\usepackage[ruled,vlined,linesnumbered]{algorithm2e}
\usepackage{subcaption}	
\usepackage{verbatim}
\usepackage{mathrsfs}
\usepackage{pgfplots}
\usepackage{dsfont}
\usepackage{siunitx}
\usepgfplotslibrary{fillbetween}

\usepackage{amsfonts,amsmath}

\newcommand{\R}{\mathbb{R}}
\newcommand{\C}{\mathbb{C}}

\usepackage{stackengine}
\def\subrangle#1{\stackengine{5pt}{}{$\!\scriptstyle #1$}{U}{l}{F}{F}{L}}
\makeatletter
\let\save@rangle\rangle
\def\rangle{\save@rangle\@ifnextchar_{\expandafter\subrangle\@gobble}{}}
\makeatother
\usepackage{nicefrac, xfrac}

\usepackage{amsthm}
\newtheorem{theorem}{Theorem}[section]
\newtheorem{definition}[theorem]{Definition}

\newtheorem{example}[theorem]{Example}

\newtheorem{remark}[theorem]{Remark}
\newtheorem{lemma}[theorem]{Lemma}
\newtheorem{problem}[theorem]{Problem}

\usepackage{array,calc}
\usepackage{tkz-euclide}
\usepackage{standalone}
\usepackage{float}
\usepackage{tikz,pgfplots}
\usepackage{tikz-cd}
\usetikzlibrary{decorations,patterns}
\usetikzlibrary{shadings}
\usetikzlibrary{fit}
\usetikzlibrary{matrix} 

\title{QR-based Parallel Set-Valued Approximation with Rational Functions}
\author{
  Simon Dirckx\\
  \texttt{simon.dirckx@austin.utexas.edu}
  \and
  Karl Meerbergen\\
  \texttt{karl.meerbergen@kuleuven.be}
  \and
  Daan Huybrechs\\
  \texttt{daan.huybrechs@kuleuven.be}
}
\date{November 2023}

\begin{document}

\maketitle

\begin{abstract}
In this article a fast and parallelizable algorithm for rational approximation is presented. The method, called (P)QR-AAA, is a (parallel) set-valued variant of the AAA algorithm for scalar functions. It builds on the set-valued AAA framework introduced by Lietaert, Meerbergen, P{\'e}rez and Vandereycken, accelerating it by using an approximate orthogonal basis obtained from a truncated QR decomposition. We demonstrate both theoretically and numerically this method's accuracy and efficiency. We show how it can be parallelized while maintaining the desired accuracy, with minimal communication cost.
\end{abstract}

\section{Introduction}\label{sec:intro}
Set-valued rational approximation, also called vector-valued rational approximation, is a recent, powerful tool for the approximation of matrix-valued functions and large sets of functions simultanueously (see, e.g., \cite{HochmanFastAAA},\cite{SuBaiSVRat}, \cite{robustRatGuttel} and \cite{KarlSVAAA}). The underlying approach to rational approximation is typically a greedy iterative scheme called `AAA' (\emph{Adaptive Antoulas-Anderson}) representing the rational functions in barycentric form (see \cite{scalarRationalAAAAntoulas} and \cite{AAA}). The central idea in set-valued rational approximation is to maximally exploit structure present in a set-valued function $z \mapsto \mathbf{f}(z)=[f_1(z),\ldots,f_N(z)]$, by constructing an approximation $\mathbf{r}(z)=[r_1(z),\ldots,r_{N}(z)]$ such that each component $r_i(z),i\in\{1,\ldots,N\}$ is a rational function with the same poles and interpolation points as every other component. We will make this idea more precise in Section \ref{sec:introSVAAA}. Typically, $\mathbf{r}(z)$ has a larger rational degree than any individual rational approximation to a component $f_i$ would have had. However, the total information needed to describe $\mathbf{r}(z)$ and the cost to evaluate it are often much smaller compared to using a union of component-wise approximations. As an additional benefit, the set-valued approach allows one to weight component functions of $\mathbf{f}$ to contribute more or less to a pre-determined cost function for approximation. For example, in \cite{robustRatGuttel}, the authors introduce \emph{weighted AAA} for matrix-valued functions in \emph{split form} (see also Section~\ref{sec:qrAAAtest}). Here, the scalar functions appearing in the split form are strategically scaled by the norms of the corresponding matrices, improving the `robustness' (as defined in \cite{robustRatGuttel}) of the final approximation. More generally, depending on analytic knowledge or predetermined design criteria, certain component functions in $\mathbf{f}$ may contribute more or less to a desired final approximation condition. Set-valued approximation in this case allows one to exploit the relative importance of the component functions, in a way that is not (easily) accomplished by approximating each function individually. Finally, the components having common interpolation points may simplify further algebraic manipulations of the approximation, such as methods for linearization and root-finding, which play an enabling role e.g. in Compact Rational Krylov (CORK) methods for nonlinear eigenvalue problems~\cite{CORK}.\\
In Section \ref{sec:introSVAAA}, we introduce the AAA algorithm and its set-valued variant, SV-AAA. Additionally we provide an interpretation of SV-AAA as an interpolative matrix decomposition. In Section \ref{sec:qr-aaa} we introduce QR-based set-valued AAA (QR-AAA) and show that it is as accurate as regular SV-AAA. In Section \ref{sec:parallel} we outline and analyze a parallel variant of QR-AAA, called PQR-AAA. We apply QR-AAA to some well-known reference non-linear eigenvalue problems in Section \ref{sec:qrAAAtest}. We conclude Section \ref{sec:qrAAAtest} by showing the effectiveness of PQR-AAA for the problem of Boundary Element Method near-field wavenumber dependence approximation (see, e.g., \cite{Dirckx}, where the far-field case was examined).

\section{Set-Valued AAA}\label{sec:introSVAAA}
In this section we outline an interpolatory rational approximation scheme known as AAA\footnote{pronounced `triple A'} (see \cite{AAA}, based on \cite{scalarRationalAAAAntoulas}), as well as its set-valued extension, called \emph{set-valued AAA} (SV-AAA, see \cite{KarlSVAAA}).\\
The AAA method is a greedy adaptive (scalar) rational approximation scheme. It is practical, easy-to-use and robust, i.e. reliable for a user-defined precision up to at least $10^{{-8}}$. Given a sufficiently fine discretization $Z$ for the approximation domain of a function $f$, in practice AAA can be treated as a black-box method. Its stability and well-conditioning rest on the \emph{barycentric form} of the rational approximant. The degree $(m-1,m-1)$ rational interpolant obtained after $m$ iterations of the AAA algorithm is given in the form
$$f(z) \approx r_m(z)=n_{m}(z)\Big/d_m(z)=\left(\sum_{\nu=1}^m\frac{w_{\nu}f(z_{\nu})}{z-z_{\nu}}\middle/\sum_{\nu=1}^m\frac{w_{\nu}}{z-z_{\nu}}\right)$$
with pairwise distinct support points $Z_m=\{z_1,\ldots,z_m\}\subset\C$ and nonzero complex barycentric weights $w_k$, subject to $\sum_\nu |w_{\nu}|^2=1$. Note that, with
\begin{equation}\label{eq:nodePoly}
\ell_m(z):=\prod_{\nu=1}^m(z-z_\nu)
\end{equation}
the \emph{node polynomial} of $Z_m$, writing $n_m/d_m=(\ell_mn_m)/(\ell_md_m)=p_{m-1}/q_{m-1}$ shows that $r_m$ is indeed a degree $(m-1,m-1)$ rational function. The expression for $r_m$ is understood to be taken in the `H\^opital limit' as $z\to Z_m$.

The greedy aspect of the AAA algorithm lies in its choice of the support point $z_m$ at iteration $m$. Setting $Z\subset\C$ as the discrete set of possible support points, the $m$th AAA approximant $r_m$ to $f$ is defined by
\begin{itemize}
\item $z_{m}:=\arg\max_{Z\backslash Z_{m-1}} \left|f(z)-r_{m-1}(z)\right|$,
\item $\|d_{m}f-n_{m}\|_{2}$ is minimal over $Z \backslash Z_{m}$.
\end{itemize}
Here, we again impose $\sum_\nu |w_{\nu}|^2=1$. To ensure the minimization property, a least-squares system must be solved at each iteration. Setting $Z\backslash Z_{m}=\{\zeta_1,\ldots,\zeta_{|Z|-m}\}$, this is achieved by computing the compact singular value decomposition (SVD) of the Loewner matrix
$$L^{(m)}:=\begin{bmatrix}
\frac{f(\zeta_1)-f(z_1)}{\zeta_1-z_1}&\cdots&\frac{f(\zeta_1)-f(z_m)}{\zeta_1-z_m}\\
\vdots&\ddots&\vdots\\
\frac{f(\zeta_{|Z|-m})-f(z_1)}{\zeta_{|Z|-m}-z_1}&\cdots&\frac{f(\zeta_{|Z|-m})-f(z_m)}{\zeta_{|Z|-m}-z_m}
\end{bmatrix}\in\C^{(|Z|-m)\times m}$$
and subsequently setting the weight vector $\mathbf{w}:=[w_1,\ldots,w_m]^T$ to be the singular vector associated with the smallest (not necessarily unique) singular value of $L^{(m)}$. This weight update step can be performed in $\mathcal{O}(m^2(|Z|-m))$ flops, for a total of $\mathcal{O}(m^3 |Z|)$ flops over AAAs run\footnote{Since $m$ is typically very small, the total cost associated with the weight update step is rather small.} up to iteration $m$. The AAA algorithm provides a remarkably robust, flexible and stable framework for rational approximation (see, e.g., \cite{AAAequi}, \cite{robustRatGuttel}, \cite{KarlSVAAA} and \cite{GuideSaadRatApprox}).\\
Taking this idea one step further, this paper concerns the approximation of set-valued functions
$$\mathbf{f}:\mathbb{C}\to\mathbb{C}^N:z\mapsto \mathbf{f}(z)$$
by set-valued rational functions of type $(m-1,m-1)$ given in barycentric form
\begin{equation}\label{eq:sv-aaa-approx}
\forall z\in \mathbb{C}: \mathbf{f}(z)\approx \mathbf{r}_m(z)=\mathbf{n}_{m}(z)\Big/\mathbf{d}_m(z)=\left(\sum_{\nu=1}^m\frac{w_{\nu}\mathbf{f}(z_{\nu})}{z-z_{\nu}}\middle/\sum_{\nu=1}^m\frac{w_{\nu}}{z-z_{\nu}}\right),  
\end{equation}
as in \cite{KarlSVAAA}. The components of $\mathbf{f}$ will be written as $\mathbf{f}(z):=[f_1(z),\ldots,f_N(z)]$, and similarly we write $\mathbf{r}_{m}(z)=[r_{m,1}(z),\ldots,r_{m,N}(z)]$. As before we denote by $Z_m\subset Z$ the set of chosen support points $\{z_{\nu}\}_{\nu=1}^m$. Equation~\eqref{eq:sv-aaa-approx} implies that the final degree $(m-1,m-1)$ rational approximations $r_{m,1},\ldots,r_{m,N}$ to $f_1,\ldots,f_N$ share support points and weights, and hence they also share poles.
\\
As in the case of regular AAA, the SV-AAA algorithm takes as an input a routine to evaluate $\mathbf{f}$ and a candidate set $Z:=\{z_1,\ldots,z_{|Z|}\}\subset \C$, and greedily selects at every iteration $m$ a new support point $z_m$ by minimizing
$$z_{m}:=\arg\max_{Z\backslash Z_{m-1}} \|\mathbf f(z)-\mathbf{r}_{m-1}(z)\|_{p}.$$
for some $\|\cdot\|_{p}$-norm. Then, it determines the corresponding weights that minimize the linearized residual $\|\mathbf{d}_{m}\mathbf{f}-\mathbf{n}_{m}\|_{p}$ over $Z \backslash Z_{m}$. The two natural choices of the $\|\cdot\|_{p}$-norm are $p=2$ and $p=\infty$. To find the weights, again a singular vector associated with the smallest singular value of a block Loewner matrix must be computed, which is of the form

$$L^{(m)}:=\begin{bmatrix}
\frac{f_1(\zeta_1)-f_1(z_1)}{\zeta_1-z_1}&\cdots&\frac{f_1(\zeta_1)-f_1(z_m)}{\zeta_1-z_m}\\
\vdots&\ddots&\vdots\\
\frac{f_1(\zeta_{|Z|-m})-f_1(z_1)}{\zeta_{|Z|-m}-z_1}&\cdots&\frac{f_1(\zeta_{|Z|-m})-f_1(z_m)}{\zeta_{|Z|-m}-z_m}\\
&&&\\
\frac{f_2(\zeta_1)-f_2(z_1)}{\zeta_1-z_1}&\cdots&\frac{f_2(\zeta_1)-f_2(z_m)}{\zeta_1-z_m}\\
\vdots&\ddots&\vdots\\
\frac{f_2(\zeta_{|Z|-m})-f_2(z_1)}{\zeta_{|Z|-m}-z_1}&\cdots&\frac{f_2(\zeta_{|Z|-m})-f_2(z_m)}{\zeta_{|Z|-m}-z_m}\\
&&&\\
\vdots&\vdots&\vdots\\
&&&\\
\frac{f_N(\zeta_1)-f_N(z_1)}{\zeta_1-z_1}&\cdots&\frac{f_N(\zeta_1)-f_N(z_m)}{\zeta_1-z_m}\\
\vdots&\ddots&\vdots\\
\frac{f_N(\zeta_{|Z|-m})-f_N(z_1)}{\zeta_{|Z|-m}-z_1}&\cdots&\frac{f_N(\zeta_{|Z|-m})-f_N(z_m)}{\zeta_{|Z|-m}-z_m}
\end{bmatrix}\in\C^{N(|Z|-m)\times m}.
$$
Computing the SVD of $L^{(m)}$ naively would result in an $\mathcal{O}(m^2N(|Z|-m))$ cost at every iteration. However, in \cite{KarlSVAAA}, a strategic updated SVD scheme for the Loewner matrix is used, reducing the cost for the weight update step at each AAA iteration to $\mathcal{O}(m^2N+m^3)$ FLOP. Since in typical applications $m\ll |Z|$, this can greatly reduce the overal computation time. The implementation in \cite{KarlSVAAA} of SV-AAA results in a cost at iteration $m$ of $\mathcal{O}(m^3+m^2N+m|Z|N)$, with modest prefactors.
This shows that the main disadvantage of set-valued AAA is that very large collections of functions to be approximated result in a prohibitively high computational cost. This is for instance the case whenever $\mathbf{f}$ represents a matrix-valued function (e.g., in wave-number dependent BEM matrices, nonlinear eigenvalue problems,\ldots) or in big-data applications (e.g., RGB image processing, sensor data analysis, digital twins,\ldots).
\subsection{SV-AAA as a Row Interpolative Decomposition}
We close this section with a useful reframing of SV-AAA approximation as a matrix approximation problem. Suppose a set-valued function $\mathbf{f}(z)=[f_1(z),\ldots,f_N(z)]$ is given, and a set-valued AAA approximant $\mathbf{r}_m(z)=[r_{m,1}(z),\ldots,r_{m,N}(z)]$ on $Z=\{z_1,\ldots,z_{|Z|}\}\subset\C$ is constructed. We introduce now two matrices, $F(Z)$ and $\widetilde{F}_m(Z)$ by collecting the discretized component functions of $\mathbf{f}$ and $\mathbf{r}_m$ as columns in an array. Concretely:
\begin{equation}\label{eq:defFmat}
    F(Z):=\begin{bmatrix}
    f_1(z_1)&f_2(z_1)&\cdots&f_N(z_1)\\
    f_1(z_2)&f_2(z_2)&\cdots&f_N(z_2)\\
    \vdots&\vdots&\ddots&\vdots\\
    f_1(z_{|Z|})&f_2(z_{|Z|})&\cdots&f_N(z_{|Z|})
\end{bmatrix}\in\C^{Z\times N}
\end{equation}
and
\begin{equation}\label{eq:defFmmat}
    \widetilde{F}_m(Z):=\begin{bmatrix}
    r_{m,1}(z_1)&r_{m,2}(z_1)&\cdots&r_{m,N}(z_1)\\
    r_{m,1}(z_2)&r_{m,2}(z_2)&\cdots&r_{m,N}(z_2)\\
    \vdots&\vdots&\ddots&\vdots\\
    r_{m,1}(z_{|Z|})&r_{m,2}(z_{|Z|})&\cdots&r_{m,N}(z_{|Z|})
\end{bmatrix}\in\C^{Z\times N}.
\end{equation}
When the set $Z$ is clear from context or not important, we will abbreviate $F:=F(Z)$ and $\widetilde{F}_m:=\widetilde{F}_m(Z)$. While selecting a column of $F$ (or $\widetilde{F}_m$) corresponds to evaluating a component function on $Z$, selecting a row can be interpreted as evaluating the set-valued function $\mathbf{f}$ (or $\mathbf{r}_m$) at a single point in $Z$, i.e. $F(i,:)=\mathbf{f}(z_i)^T$. For this reason, we will use the terminology `SV-AAA approximation of $\mathbf{f}$' and `SV-AAA approximation of $F$' interchangeably. We will also make the following assumption:
$$F\text{ is scaled such that }\forall j: \|F(:,j)\|_{\infty}=1$$
which can always be achieved by applying a diagonal scaling to $F$, and possibly discarding zero columns. This way, the error results in this section and in sections~\ref{sec:qr-aaa} and~\ref{sec:parallel} are automatically relative error results.\\
An approximate \emph{row interpolative decomposition} (RID) of a matrix $A\in\C^{n_1\times n_2}$, with respect to some matrix norm $\|\cdot\|$ is a decomposition of the form
$$A\approx A_m=H_m\cdot A(I_m,:)$$
with $H_m\in\C^{n_1\times m}$ and $I_m\subseteq\{1,\ldots,q\}$, such that $\|A-A_m\|$ is sufficiently small. These have been studied in, e.g., \cite{DongID}. The first important observation of our manuscript, captured in Theorem~\ref{thm:obs1}, is that an SV-AAA approximation can be interpreted as an RID of the matrix $F$ from equation~\eqref{eq:defFmat}.
\begin{theorem}\label{thm:obs1}Suppose $\mathbf{f}:\C\to\C^N$ can be approximated on $Z\subset\C$ by the set-valued rational function $\mathbf{r}_m$ with support points $Z_m \subset Z$ as in Equation~\eqref{eq:sv-aaa-approx}, such that for some $\epsilon>0$ 
\begin{equation}\label{eq:eqobs1}
\textbf{res}_m:=\sup_{z\in Z}\|\mathbf{f}(z)-\mathbf{r}_m(z)\|_{p}<\epsilon.
\end{equation}
Let $F:=F(Z)$ and $\widetilde{F}_m:=\widetilde{F}_m(Z)$ be as in equations~\eqref{eq:defFmat} and~\eqref{eq:defFmmat}. Then the matrix $\widetilde{F}_m$ is of rank $m$, and can be written as
$$\widetilde{F}_m=H_m\cdot F(Z_m,:).$$
Thus, $F\approx \widetilde{F}_{m}=H_m\cdot F(Z_m,:)$ constitutes an approximate RID with respect to the \emph{row-wise} $\|\cdot\|_{p,\infty}$-norm, by which we mean
\begin{equation}\label{eq:normStatementRID}
    \|F-\widetilde{F}_m\|_{p,\infty}:=\max_{i\in |Z|}\|F(i,:)-\widetilde{F}_m(i,:)\|_p<\epsilon.
\end{equation}
\end{theorem}
\begin{proof}
Recall from equation~\eqref{eq:sv-aaa-approx} that $$\mathbf{r}_m(z) = \left(\sum_{\nu=1}^m\frac{w_{\nu}\mathbf{f}(z_{\nu})}{z-z_{\nu}}\middle/\sum_{\nu=1}^m\frac{w_{\nu}}{z-z_{\nu}}\right)$$ with $\{z_{\nu}\}_{\nu}^m=Z_m$ the support points selected by SV-AAA. Using this, the matrix $\widetilde{F}_m$ can be written as
$$\widetilde{F}_m=\sum_{\nu=1}^m\mathbf{h}_{\nu}\mathbf{f}(z_{\nu})^T$$
with $\mathbf{h}_{\nu}\in\mathbb{C}^{Z}$ defined by $\mathbf{h}_{\nu}(i):=\left.\frac{w_{\nu}}{z_i-z_{\nu}}\middle/\left(\sum_{j=1}^m\frac{w_{j}}{z_i-z_{j}}\right)\right.$, again evaluated as a limit for $z_i\in Z_m$. Setting $H_m = [\mathbf{h}_{1},\ldots,\mathbf{h}_{m}]\in\C^{Z\times m}$ we can re-write this as
$$\widetilde{F}_m = H_m\cdot F(Z_m,:).$$Equation~\eqref{eq:normStatementRID} follows from the requirement in equation~\eqref{eq:eqobs1} and the definition of the $\|\cdot\|_{p,\infty}$-norm.
\end{proof}
The actual error criterion we are interested in for SV-AAA approximation of $F$ is not necessarily $\|F-\widetilde{F}_m\|_{p,\infty}<\epsilon$, but rather the more useful $\|F-\widetilde{F}_m\|_{\max}<\epsilon$ implied by it, because this means all functions in $\mathbf{f}$ can be expected to be well-approximated in the uniform norm.\\
The final key insight that results in our QR-AAA method, captured in Theorem~\ref{thm:transitive}, is what we call the \emph{transitive property} of RIDs. Together with Remark~\ref{rem:TransativeProp}, Theorem~\ref{thm:transitive} is a generalization of results in \cite{VoroninMartinsson} and \cite{DongID}.
\begin{theorem}[transitive property of RIDs]\label{thm:transitive}
    Let $A\in \C^{n_1\times n_2}$ be a matrix that can be factorized as $A=XY$, $X\in\C^{n_1 \times k}$ and $Y\in\C^{k \times n_2}$ and let $X\approx X_m=H_m\cdot X(I_m,:)$ be an approximate RID for $X$ w.r.t. the norm $\|\cdot\|_{p,\infty}$ such that
    $$\|X-X_m\|_{p,\infty}<\epsilon.$$
    Then $A\approx A_m = H_m\cdot A(I_m,:)$ is an approximate RID for $A$ w.r.t. $\|\cdot\|_{\max}$ such that
    $$\|A-A_m\|_{_{\max}}<\frac{k\epsilon}{\sqrt[\uproot{3}p]{k}}\|Y\|_{\max}$$
\end{theorem}
\begin{proof}
    As $X(I_m,:)Y=A_m(I_m,:)$, we have $X_mY=H_mA(I_m,:)$. Additionally, $\|A-A_m\|_{\max}=\|(X-X_m)Y\|_{\max}\leq k^{1-1/p}\|X-X_m\|_{p,\infty}\|Y\|_{\max}<\epsilon k^{1-1/p}\|Y\|_{\max}$, since for any two vectors $\mathbf{x},\mathbf{y}\in\C^k$ we have (by H\"older's inequality) 
    $$|\mathbf{x}^T\mathbf{y}|\leq\|\mathbf{x}\|_1\|\mathbf{y}\|_{\infty}\leq k^{1-1/p}\|\mathbf{x}\|_{p}\|\mathbf{y}\|_{\infty},$$
    which concludes the proof.
\end{proof}
Because in Theorem~\ref{thm:transitive} the matrix $H_m$ is the same for the RID of $X$ and $A$, a low-rank decomposition $F\approx XY$ (with $F$ as in equation~\eqref{eq:defFmat}), followed by an SV-AAA approximation on $X$, actually constitutes a set-valued rational approximation of $F$ which is interpolative precisely because of the transitive property.\\
This is the basic principle behind QR-AAA, formalized in section~\ref{sec:qr-aaa}, where $X=Q$ and $Y=R$ correspond to a rank-revealing QR-decomposition of $F$.
\begin{remark}\label{rem:TransativeProp}
    Theorem~\ref{thm:transitive} holds more generally. Suppose $\|\cdot\|$, $\|\cdot\|_{\alpha}$ and $\|\cdot\|_{\beta}$ are matrix norms such that for all $M_1\in\C^{n_1\times k}$, $M_2\in\C^{k\times n_2}$ there is a constant $c_{\alpha,\beta}(k)$ independent of $M_1,M_2$ for which $\|M_1M_2\|\leq c_{\alpha,\beta}(k)\|M_1\|_{\alpha}\|M_2\|_{\beta}$. Then, with $\|X-X_m\|_{\alpha}<\epsilon$, the conclusion of Theorem~\ref{thm:transitive} still holds, but with $\|A-A_m\|<\epsilon c_{\alpha,\beta}(k)\|Y\|_{\beta}$.
\end{remark}

\section{QR-based SV-AAA}\label{sec:qr-aaa}
In this section we outline `QR-based SV-AAA' (QR-AAA). We will show that, under suitable assumptions, QR-AAA guarantees accuracy up to user-specified tolerance. For the remainder of this text we will set $\|\cdot\|_{p,\infty}=\|\cdot\|_{2,\infty}$. However, all results can be generalized to the case $p\in [1,\infty]$.\\
In short, QR-AAA is a rank-revealing QR-factorization with column pivoting (RRQR), followed by a weighted SV-AAA approximation of $Q$ (not to be confused with weighted AAA from \cite{robustRatGuttel}). In our implementation, the RRQR step is performed using Householder reflectors (see \cite[Chapter~5]{GolubAndVanLoan}). Given a matrix $A\in\mathbb{C}^{n_1\times n_2}$, RRQR returns factors $Q\in\C^{n_1 \times k}$ and $R\in\C^{k\times n_2}$ and a permutation $\Pi$ of $\{1,\ldots,n\}$ such that
$A\Pi\approx QR$ up to user-specified tolerance, $Q^*Q=\mathrm{1}_k$ and $R$ is upper triangular.
In addition, we have the following properties, which we include without proof:
\begin{lemma}\label{lem:qr-properties}
Suppose $A\Pi\approx QR$, $Q\in\C^{n_1\times k}$, $R\in\C^{k \times n_2}$ is the output of RRQR with column pivoting for some tolerance $\epsilon$. Then
$$|R(1,1)|\geq|R(2,2)|\geq\cdots\geq|R(k,k)|>0$$
and
$$|R(i,i)|>|R(i,j)|$$
for $1\leq i\leq k$ and $i\leq j\leq N.$ Additionally, $\|A\Pi-QR\|_{\max}<\epsilon$ and the computational cost of RRQR is $\mathcal{O}(n_1n_2k+k^3)$ (see \cite[Chapter~5]{GolubAndVanLoan}).
\end{lemma}
At a tolerance $\epsilon$, the matrix $Q=[\mathbf{q}_1,\ldots,\mathbf{q}_k]$ returned by RRQR represents the so-called \emph{$\epsilon$-numerical range} of $A$. The diagonal elements $R(i,i)$ of $R$, also called the \emph{pivots}, intuitively indicate the importance of each $\mathbf{q}_i$ in the linear mapping defined by $A$.

\begin{algorithm}[ht]
\SetAlgoLined
\SetKwInOut{Input}{input}
\SetKwInOut{Output}{output}
\SetKwInOut{Init}{init}
\caption{QR-based SV-AAA (QR-AAA)}\label{alg:qr-aaa}
\Input{$F\in \mathbb{C}^{Z \times N}$ (without zero columns), convergence tolerance $\epsilon$}
\Output{$(Z_m,W_m) \in Z^{m}\times\mathbb{C}^{m}$, supports and weights of SV-AAA, such that $\|F-H_mF(Z_m,:)\|_{\text{max}}<\epsilon$, with $H_m$ as in Theorem~\ref{thm:obs1}}
\Init{$G:=0\in \mathbb{C}^{Z\times N}$}
\vspace{.5em}
\tcc{Part 1: scale F}
\For{$j=1,\ldots,N$}{
    $G(:,j) \leftarrow F(:,j)/\|F(:,j)\|_{\infty}$\;
}
\tcc{Part 2: Construct $Q$ and $R$}
$[Q,R]\leftarrow \textbf{RRQR}(G,\epsilon/2)$\;
$k\leftarrow \text{size}(Q,2)$\;
\tcc{Part 3: Construct $(\widetilde{Q}_m,Z_m,W_m)$}
$\Gamma \leftarrow diag(R)$;\\
$(\widetilde{Q}_m,Z_m,W_m)\leftarrow \textbf{SV-AAA}(Q\Gamma,\epsilon/2\sqrt{k})$;
\end{algorithm}
Lemma~\ref{lem:qr-properties} provides us with the final insight needed for QR-AAA. If $\Gamma:=diag(R)$, then $\|\Gamma^{-1}R\|_{\max}=1$. Let $F\Pi\approx QR$ be the rank $k$ output of RRQR with column pivoting, i.e. $\|F\Pi - (Q\Gamma)(\Gamma^{-1}R)\|_{\max}<\epsilon$, and $$Q\Gamma \approx H_m\cdot(Q\Gamma)(Z_m,:)=(H_m\cdot Q_m(Z_m,:))\Gamma=:\widetilde{Q}_m\Gamma$$
be obtained by applying SV-AAA to $Q\Gamma$, with tolerance $\epsilon$. By Theorem~\ref{thm:transitive}
\begin{equation}\label{eq:errEQ}
    \|F\Pi-\widetilde{Q}_mR\|_{\max}\leq \epsilon + \sqrt{k}\|Q\Gamma-\widetilde{Q}_m\Gamma\|_{2,\infty}\|\Gamma^{-1}R\|_{\max}<\epsilon(1+\sqrt{k}).
\end{equation}
The process of approximating $Q\Gamma$ rather than $Q$ is what we mean by `weighted SV-AAA' applied to $Q$. As shown in section~\ref{sec:introSVAAA}, by Theorem~\ref{thm:obs1} the support nodes and weights $(Z_m,W_m)$ returned by SV-AAA applied on $Q\Gamma$ can now be used as support points and weights for a set-valued rational approximation of $F\Pi$, and hence also of $F$.

\begin{figure}
\centering
\includegraphics[width=.8\linewidth]{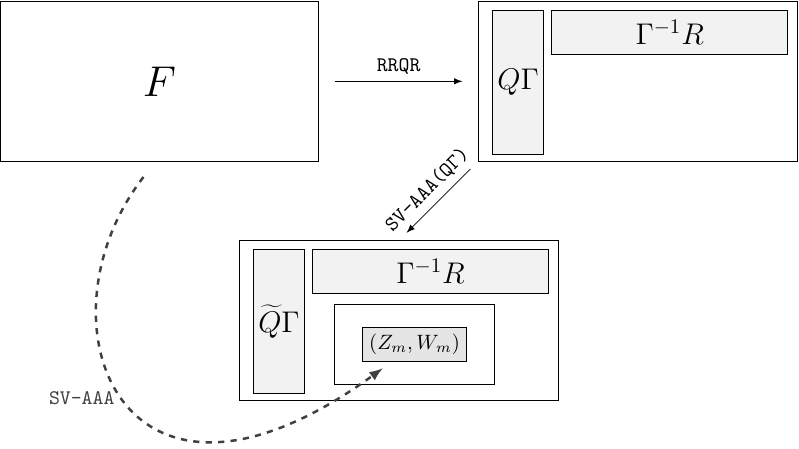}
\caption{Diagram showing the principle of QR-AAA. Here $\Gamma=\text{diag}(R)$ and $\widetilde{Q}=H_mQ(Z_m,:)$ with $H_m$ as in Theorem~\ref{thm:obs1}. The dashed \texttt{SV-AAA} arrow indicated that the final $(Z_m,W_m)$ constitute a set-valued rational approximation $F\approx H_m F(Z_m,:)$.}\label{fig:qrAAA}
\end{figure}

The implementation of QR-AAA, which we take as its definition, is given in Algorithm~\ref{alg:qr-aaa}. The tolerances for the RRQR step and the SV-AAA step are reduced to $\epsilon/2$ and $\epsilon/2\sqrt{k}$ respectively, to ensure an overall error of $\epsilon$ (see equation~\eqref{eq:errEQ}).
Figure \ref{fig:qrAAA} shows a simple diagram illustrating QR-AAA. Note that the final approximation $F\approx H_mF(Z_m,:)$ is not identical to the output of SV-AAA applied to $F$, but constitutes a different set-valued rational interpolatory approximation, in barycentric form, that well-approximates $F$. It can therefore be said to behave like a SV-AAA approximation of $F$, hence the dashed arrow labeled `SV-AAA'.\\
Using the diagonally scaled $Q\Gamma$ is essential, as the trailing columns of $Q$ tend to become more oscillatory for increased precision $\epsilon$, leading to an increased rational degree in the SV-AAA approximation of $Q$ if unweighted by $\Gamma$.
In Section \ref{sec:qrAAAtest} we report the remarkable speed-ups attained by QR-AAA.
\begin{remark}\label{rem:QRexists}
The idea of using a QR decomposition in approximation algorithms is not new. Typically, the orthogonal bases obtained by the QR decomposition provide additional stability, as in the RBF-QR method (see, e.g.,~\cite{RBFQR} and~\cite{RBFflat}). In fact, it is standard practice to orthogonalize bases of discretized functions using, e.g., Gram-Schmidt or the QR decomposition. Nevertheless, our results on the accuracy of set-valued approximation using the truncated QR decomposition seem new, as is the observation of its computational benefits in the context of set-valued approximation.
\end{remark}
\begin{remark}\label{rem:diffBasis}
If from analytical knowledge, some desirable basis $\{\varphi_j\}_{j=1}^k$ whose linear span (approximately) contains $\{f_i\}_{i=1}^N$ is known, one can replace $Q$ above by $\Phi\in\C^{Z\times k}$, with $\Phi(i,j)=\varphi_j(z_i)$. We can then write $F\approx \Phi C$, where $C=\Phi^{\dagger}F$. In this case, the principle of QR-AAA still works, with $Q$ replaced by $\Phi$ and $\Gamma$ the diagonal matrix defined by $\Gamma_{ii}:=\|C(i,:)\|_{\infty}$.
\end{remark}
A \texttt{C++} implementation of QR-AAA has been made publicly available at \cite{PQRAAAGit}.

\section{Parallel QR-AAA}\label{sec:parallel}
While QR-AAA already performs significantly better than SV-AAA for many problems, it is conceivable that for extremely large problems, i.e., when $F$ does not fit in memory, the QR-AAA method can still be prohibitively costly. Additionally, the dominant cost in SV-AAA and QR-AAA is often the computation of $F$. For these reasons, we can look for methods to exploit today's increasing parallel computing power. Below, one such method for parallel QR-AAA, which we call PQR-AAA\footnote{In keeping with the alphabet theme} is outlined. The proposed approach is as follows. We start by distributing (the computation of) $F=[F_1,\ldots,F_{l}]$ over $l$ processors\footnote{We will restrict ourselves to the case where $F_{\mu}$ fits info the memory of node $\mu$. The methods described in this section generalize to parallel nodes running multiple threads.}. The first two steps for each processor $\mu$ are: 
\begin{itemize}
    \item[(1)] Processor $\mu$ computes $F_{\mu}\approx Q_{\mu}R_{\mu}$.
    \item[(2)] Processor $\mu$ applies QR-AAA to construct $\widetilde{Q}_{\mu}$, by computing supports and weights $(Z_{\mu},W_{\mu})$.
\end{itemize}
Here, the subscript `$\mu$' refers to the processor index, not the rational degree. Without loss of generality we will assume for the remainder of this section that the degree of rational approximation of each $\widetilde{Q}_{\mu}$ is $m-1$.
\\
After steps $(1)$ and $(2)$ we have a piece-wise set-valued rational approximant, consisting of $[\widetilde{Q}_1\cdots\widetilde{Q}_l]$. The next goal is to `glue' these approximations together, which we call the \texttt{accumulate} step described in more detail further on. Of course, we do not want to re-approximate $[\widetilde{Q}_1\cdots\widetilde{Q}_l]$ over all of $Z$ again -- this defeats the purpose of parallel approximation. We are thus faced with the following problem:
\begin{problem}[Accumulate support problem]
    Given a set of set-valued rational functions $\{\widetilde{Q}_{\mu}\}_{\mu}$ with support nodes and weights $\{Z_{\mu},W_{\mu}\}$ such that $Z_{\mu} \subset Z$, find a global set-valued rational approximation $\widetilde{Q}$ using only samples in $Z^+\subset Z$, where $Z^+$ is as small as possible.
\end{problem}

\subsection{The accumulate support problem}

The selection of $Z^+$ is crucial. It should be small enough to realize computational gain compared to using all of $Z$, but large enough so as not to introduce new errors in the \texttt{accumulate} step. Naively, one might minimally try to take $Z^+=\cup_{\mu} Z_{\mu}$ as the union of the existing support sets. However, for rational interpolation of degree $m-1$, at least $2m-1$ sample points are required.\footnote{This can be seen by counting degrees of freedom in the barycentric form: in addition to $m$ support points, $m$ weights are needed, normalized to have unit norm.} Since $|\cup_{\mu} Z_{\mu}|\geq 2m-1$ often does not hold, the set has to be extended by some $Z^{\textnormal{e}}\subset Z$. A larger set that ensures uniqueness of the interpolant may still not be sufficient,  because approximating a degree $m-1$ rational function on $\cup_{\mu} Z_{\mu}\cup Z^e$ may be very ill-conditioned. Indeed, given some rational function $f$ on $Z$, it is possible that its rational approximant $r$ on $Z^+:=(\cup_{\mu}Z_{\mu}\cup Z^{\textnormal{e}})$ satisfies 
$$\|f-r\|_{Z^{+},\infty}<\epsilon$$
while the error $\|f-r\|_{Z,\infty}$ on the full grid grows exponentially large in $m$.
\\
Currently, no theory to construct general extension sets exists. For the important case of $Z$ a finely sampled equispaced set of points in an interval $[-1,1]$, which can be adapted to $[a,b]$ or $\imath[a,b]$, we observe in our experiments that adding sufficiently many points sampled uniformly at random from $Z\,\backslash (\cup_{\mu}Z_{\mu})$ results in a satisfactory final accuracy. In essence, the points added can be viewed as a validation set for the final rational approximation.\\
It is possible to be more deliberate in the selection of $Z^{e}$. We use the theory from~\cite{AdcockOptimSampling} for bounded polynomial growth. First we fix some definitions.
\begin{definition}\label{def:finDefQRAAA}
Given a set of QR-AAA approximations $\{\hat{Q}_{\mu}\}_{\mu=1}^l$ on a common set $Z\subset\C$, with support points and weights $\{(Z_{\mu},W_{\mu})\}_{\mu=1}^{l}$, we set
$$m^+:=2\left|\cup_{\mu=1}^{l}Z_{\mu}\right|-2.$$
Furthermore, for any ordered distribution of points $Z^{\textnormal{e}}=\{z_{i}\}_{i=1}^{M}$ in $[-1,1]$, define
$$\zeta(Z^{\textnormal{e}}) := \max_i\int_{z_i}^{z_{i+1}}\frac{1}{\sqrt{1-x^2}}\dif x.$$ We say that $Z^{\textnormal{e}}\subset Z$ is a \emph{good extension set for degree $m^+$} if $m^+\zeta(Z^{\textnormal{e}})=\alpha$ with $0<\alpha<1$ sufficiently small. Additionally, we set
\begin{equation}\label{eq:sumQ}
Z^+:=\left(\cup_{\mu}Z_{\mu}\right)\cup Z^{\textnormal{e}}
\end{equation}
and define
$$Q_1\uplus\cdots\uplus Q_l:=[Q_{1}(Z^+,:)\cdots Q_{l}(Z^+,:)]$$
where the dependence of $Z^+$ on $Z^{\textnormal{e}}$ is assumed known from context. 
\end{definition}
As shown in \cite{AdcockOptimSampling}, if $\alpha<1$ we have that
\begin{equation}\label{eq:B(M,N)}
B(Z^{\textnormal{e}},m^+):=\sup_{p\in\mathbb{P}_{m^+}}\{\|p\|_{\infty,[-1,1]}\,|\,\|p\|_{\infty,Z^{\textnormal{e}}}\leq1\}
\end{equation}
is bounded by
$$B(Z^{\textnormal{e}},m^+)<\frac{1}{1-\alpha}.$$
\begin{example}\label{ex:mockCheb}
Given an equispaced grid $Z$ in $[-1,1]$, the \emph{mock Chebyshev} points (see \cite{BoydMockCheb}) of cardinality $M$, denoted by $Z_{\text{cheb},M}$ are the $M$ points in $Z$ closest to the Chebyshev nodes of order $M$. For $|Z|=\mathcal{O}(M^2)$ (see \cite{AdcockOptimSampling}, Section 5.3), the mock Chebyshev points achieve $\zeta\approx\frac{\pi}{M}$. They are useful to characterize what we mean by a `sufficiently fine' equispaced grid. We say that $Z$ is `sufficiently fine' with respect to $M$ if the mock Chebyshev points in $Z$ achieve $\zeta(Z_{\text{cheb},M})=\frac{\pi}{M}(1+\delta)$, with $0<\delta<1$ sufficently small. In that case we have, if $M=\lceil 3\pi m^+\rceil$, that
$$\frac{1}{1-\alpha}<\frac{1}{1-\frac{1+\delta}{3}}=\frac{3}{2-\delta}<3.$$
It follows that
$$B(Z_{\text{cheb},M},m^+)<3.$$
\end{example}
With equation \eqref{eq:B(M,N)} in mind we now have Theorem~\ref{thm:finalThmSVAAA}, which guarantees satisfactory approximation of the final PQR-AAA step:
\begin{itemize}
\item[(3)] each node $\mu$ extends its approximation to $Z^+ = (\cup_{\mu}Z_{\mu})\cup Z^{\textnormal{e}}$
\item[(4)] a final SV-AAA step on $\uplus_{\mu}Q_{\mu}$ is applied.
\end{itemize}
\begin{theorem}\label{thm:finalThmSVAAA}
We use the notation from Definition~\ref{def:finDefQRAAA}. Suppose $f:=n/d=p/q$ is a rational function on $[-1,1]$ of degree $m-1$, defined by support points and weights $(Z_m,W_m)$. Then if $Z_m^{+}=Z_m\cup Z^{\textnormal{e}}$, we have that any AAA approximant $\hat{f}=\hat{n}/\hat{d}=\hat{p}/\hat{q}$ of degree at most $m-1$ such that \begin{equation}\label{eq:errCritFinalTheorem}
\|\hat{d}f-\hat{n}\|_{Z_m^{+},\infty}<\epsilon
\end{equation}
satisfies
$$\|p\hat{q}-\hat{p}q\|_{\infty,[-1,1]}<\|q\|_{[1,-1],\infty}\|\hat{\ell}\|_{[1,-1],\infty}B(Z^{\textnormal{e}},m^+)\epsilon$$
where $\hat{\ell}$, as in equation~\eqref{eq:nodePoly}, is the node polynomial for $\hat{f}$.
\end{theorem}
\begin{proof}
Clearly $p\hat{q}-\hat{p}q$ is a polynomial of degree at most $m^+$. It holds that
\begin{align*}
\|p\hat{q}-\hat{p}q\|_{Z_m^+,\infty}&=\|q\hat{\ell}(f\hat{d}-\hat{n})\|_{Z_m^+,\infty}\\
&<\|q\hat{\ell}\|_{\infty,Z_m^+}\epsilon
\end{align*}
from which
$$\|p\hat{q}-\hat{p}q\|_{\infty,[-1,1]}<\|q\|_{\infty,Z_m^+}\|\hat{\ell}\|_{\infty,Z_m^+}B(Z^{\textnormal{e}},m^+)\epsilon$$
follows by definition of $B(Z^{\textnormal{e}},m^+)$ and $Z_m^+$.
\end{proof}
We make a few remarks regarding this result. First, the factor $\|q\|_{\infty,Z_m^+}<\|q\|_{\infty,[-1,1]}$ merely serves to make the error relative, since we can arbitrarily rescale the numerator and denominator of $f$ by a common factor.

The factor $\|\hat{\ell}\|_{\infty,Z^+_m}$ is more difficult to analyze. It expresses the maximal size of the node polynomial for the denominator of $\hat{f}$. In principle, it is bounded only by $2^k$, for $k-1\leq m-1$ the degree of $\hat{f}$. In practice this growth is rarely observed. In fact, it is possible to explicitly compute $\|\hat{l}\|_{\infty,Z^+_m}$ in $\mathcal{O}(|Z^+_M|)$ FLOP at each iteration and to incorporate this factor into the convergence criterion for the approximation of $f$ by $\hat{f}$.

Finally, we note that the linearized error $\|p\hat{q}-\hat{p}q\|_{\infty,[-1,1]}$ is often a good indicator of the actual error $\|f-\hat{f}\|_{\infty,[-1,1]}$ in AAA, even though this theoretically only holds if $f$ and $\hat{f}$ have no poles close to $[-1,1]$. See \cite{RatMinMax} for an exploration of this topic. A simple counter-example is when $f=(z-\imath\epsilon)^{{-}1}$ and $\hat{f}=(z-\imath\epsilon/2)^{{-}1}$, $\epsilon>0$. Here the linearized error can be made arbitrarily small, while the actual error ($=1/\epsilon$) grows arbitrarily large. For such cases, a more refined analysis is needed.

\subsection{Schematic description and communication costs}
\begin{figure}[H]
\centering
\includegraphics[width=.8\linewidth]{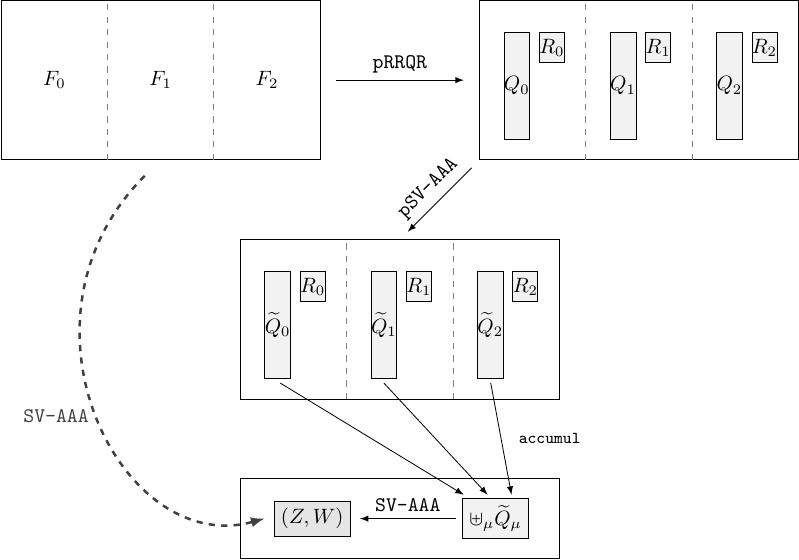}
\caption{The parallel QR based set-valued AAA approach. The original $F$ is distributed over the chosen nodes, and on each node an RRQR decomposition is executed in parallel (`\texttt{pRRQR}'). Then, in parallel, each (weighted) $Q$ is approximated using SV-AAA (`\texttt{pSV-AAA}'). Finally, these are accumulated (`\texttt{accumul}') by applying steps (3)-(4) outlined in the above.}\label{fig:paraQRSVAAA}
\end{figure}

A schematic overview of  PQR-AAA is given in Figure~\ref{fig:paraQRSVAAA}.

\begin{figure}[ht]
\centering
\scalebox{.8}{
\begin{minipage}{\linewidth}
\centering
\includegraphics[scale=.2]{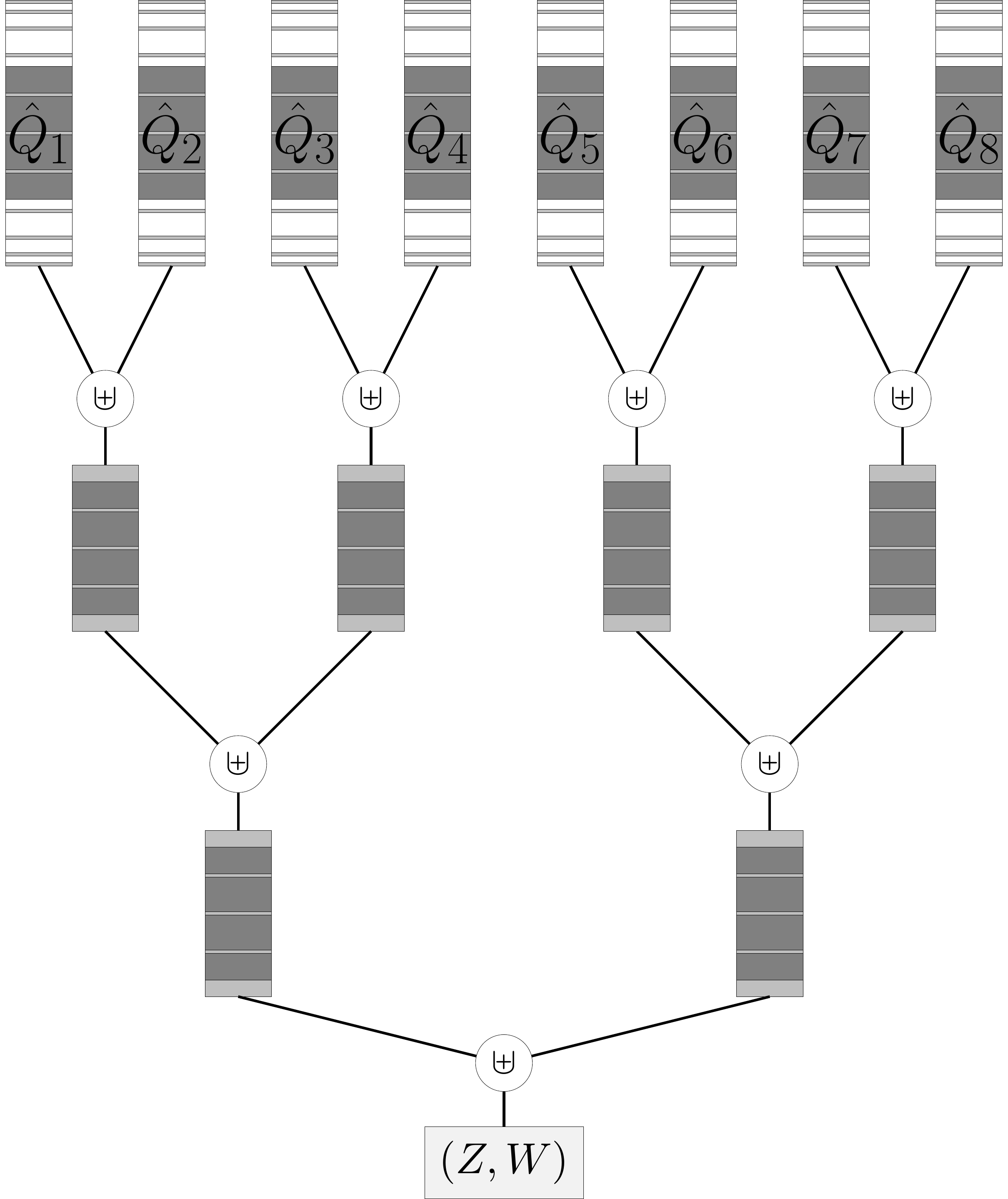}
\end{minipage}
}
\caption{Diagram illustrating the PQR-AAA pairwise accumulation stage. At each stage the `$\uplus$' sum $\hat{Q}_{\mu_1}\uplus\hat{Q}_{\mu_2}$ of the  approximations in two nodes is constructed, after which SV-AAA$(\hat{Q}_{\mu_1}\uplus\hat{Q}_{\mu_2})$ is computed.}\label{fig:accumulQRAAA}
\end{figure}

Note that PQR-AAA as explained thus far involves an all-to-all communication phase in the accumulation step. We do not know whether an all-to-all communication step for the construction of $\cup_{\mu}Z_{\mu}$, which involves the all-to-all communication of indices, can be avoided. However, importantly it is possible to avoid the more costly all-to-all communication of blocks $Q(Z_m^+,:)$. This can be achieved, e.g.,  by accumulating the computational nodes two-by-two. This is shown in Figure~\ref{fig:accumulQRAAA}.

\section{Experiments}\label{sec:qrAAAtest}
In this section we demonstrate the effectiveness of QR-AAA and PQR-AAA. We apply them first to some standard matrix functions. We close this section with an application from boundary element methods. We illustrate their convergence properties and complexity.

\subsection{NLEVP problems}
In \cite{NLEVP1} and \cite{NLEVP2}, a class of reference \emph{Non-Linear Eigenvalue Problems} (NLEVPs) was introduced, that serves as a rich source of matrix functions. The matrix functions in the NLEVP benchmark set, available at \cite{NLEVPGit} are all in \emph{split form},
i.e.
$$A(z) = \sum_{\ell=1}^L g_{\ell}(z)M_\ell$$
for some set of matrices $\{M_{\ell}\}_{\ell=1}^L$, and scalar functions $\{g_\ell\}_{\ell=1}^L$.\\
We briefly introduce the NLEVPs we use to test our method, after which we report the performance of QR-AAA. We will not go into detail about the structure and properties of the systems. More detailed information about sparsity, eigenvalue distributions, positive (semi-)definiteness of the factors etc. can be found in the NLEVP repository (see \cite{NLEVPGit}).
\subsubsection{Problem 1: A Clamped Sandwich Beam}
Our first example problem is that of the \emph{clamped sandwich beam} (\texttt{sandwich\_beam} in \cite{NLEVP2}). A full description is given in \cite{MeerbergenBeam}. The matrix function is given by $$A(s) = \left(K+\frac{G_{0}+(s\tau)^p G_{\infty}}{1+(s\tau)^p}D+s^2M\right)$$
with $s$ the Laplace variable, $p=.675$, and constants $G_0,G_{\infty},\tau\in\R$. We set $Z\subset[200\imath,30000\imath]$ equispaced of size $|Z|=1000$. We run SV-AAA and QR-AAA with tolerance $10^{\smin8}$.

\subsubsection{Problem 2: Photonic Crystals}
The second problem problem we consider originates from a 2D Helmholtz problem on photonic crystals (\texttt{photonic\_crystal} in \cite{NLEVP2}). In this case the matrix function is given by
$$A(s) = G+s^2(\epsilon_0M_0+\epsilon_1(s)M_1)$$
where $\epsilon_0\in\R$ and
$$\epsilon_1(s) = c + \sum_{\ell=1}^{L}\frac{\lambda_{P,\ell}^2}{s^2+\alpha s+\lambda_{0,\ell}^2}$$
with parameters $\alpha,c,\{\lambda_{P,\ell},\lambda_{0,\ell}\}_{\ell}$ in $\R$ and $L=2$. For more info, see \cite{photonic}. We discretize this problem with $Z\subset[0,10\imath]$ equispaced of size $|Z|=1000$. We run SV-AAA and QR-AAA with tolerance $10^{\smin8}$.

\subsubsection{Problem 3: Schr\"odinger's Equation in a Canyon Well}
Our third problem, \texttt{canyon\_particle} in \cite{NLEVP2} comes from the Schr\"odinger equation for a particle in a canyon-shaped potential well. A detailed description can be found in \cite{MeerbergenCanyon}. The resulting matrix function has the form
$$A(\lambda) = H-\lambda\mathds{1}-\sum_{\ell=1}^Le^{\imath\sqrt{m(\lambda-a_{\ell}})}S_\ell.$$
Here $\lambda$ is a real variable and $L=81$. We discretize this problem with $Z\subset[a_1+\delta,a_2-\delta]$ an equispaced of size $|Z|=1000$, where $\delta=10^{{-}4}$. We run SV-AAA and QR-AAA with tolerance $10^{\smin8}$.

\subsubsection{Problem 4: Time Delay Systems with random delays}
Our fourth problem, \texttt{time\_delay3} from \cite{NLEVPGit}, comes from the Laplace transform of a time delay system, and has the form
$$A(s)=s\mathds{1}-\sum_{\ell=1}^{L}A_{\ell}e^{-s\tau_\ell}$$
with $\{A_{\ell}\}_{\ell}$ dense random matrices representing the gain at delay $\tau_{\ell}\in\R^+$. We set $\|A_{\ell}\|_{\infty}$ to be linearly increasing with $\ell$, and we set $L=20$. Since the matrices $\{A_{\ell}\}_{\ell}$ are generated randomly, and because of the oscillatory nature of $e^{\smin s\tau_{\ell}}$, this is a tough problem for rational approximation. We set $Z\subset[-10\imath,10\imath]$ equispaced of size $|Z|=1000$. We run SV-AAA and QR-AAA with tolerance $10^{\smin4}$.\footnote{The approximation is run at a lower tolerance because for higher tolerances the memory requirements and computation time become exceedingly large.}

\subsubsection{Numerical results}

\begin{figure}[ht]
\centering
\scalebox{.8}{
\begin{minipage}{\linewidth}
\begin{subfigure}[t]{.49\linewidth}
\includegraphics[scale=.5]{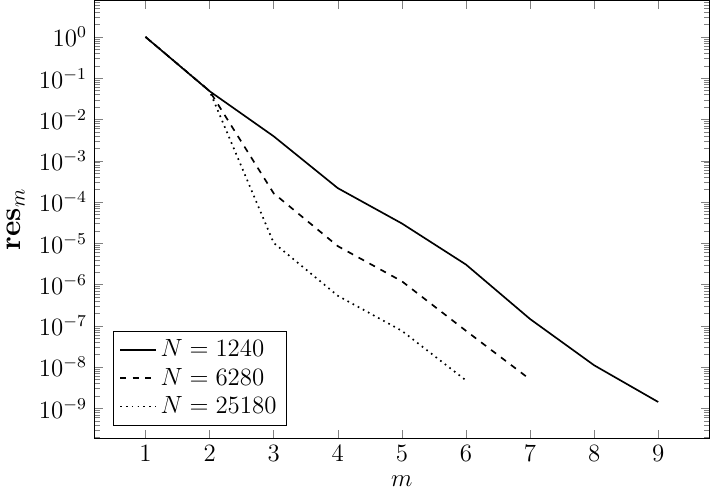}
\caption{\texttt{sandwich\_beam}: SV-AAA}
\end{subfigure}
\hfill
\begin{subfigure}[t]{.49\linewidth}
\includegraphics[scale=.5]{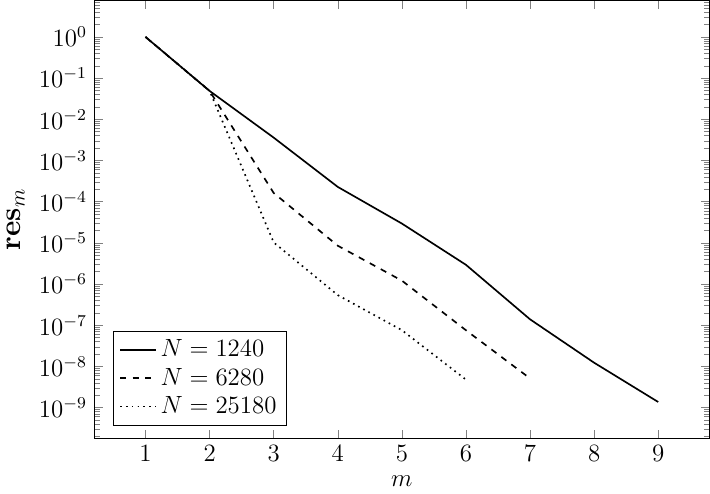}
\caption{\texttt{sandwich\_beam}: QR-AAA}
\end{subfigure}

\vspace{\baselineskip}
\begin{subfigure}[t]{.49\linewidth}
\includegraphics[scale=.5]{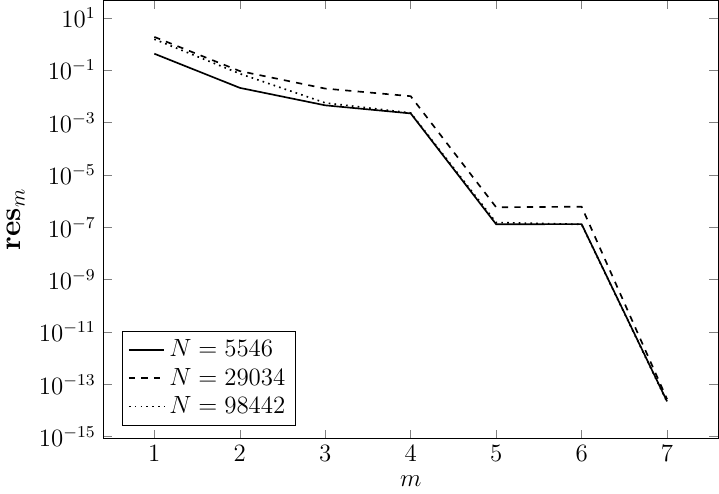}
\caption{\texttt{photonic\_crystal}: SV-AAA}
\end{subfigure}
\hfill
\begin{subfigure}[t]{.49\linewidth}
\includegraphics[scale=.5]{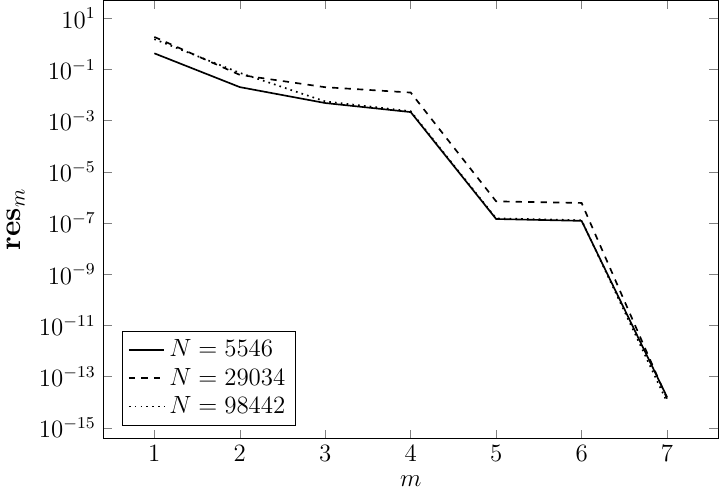}
\caption{\texttt{photonic\_crystal}: QR-AAA}
\end{subfigure}
\hfill
\end{minipage}
}
\caption{Residue $\textbf{res}_m$ for SV-AAA (left) and QR-AAA (right) over the degree $m$, for the first two selected problems.}\label{fig:errQRAAAPart1}
\end{figure}

\begin{figure}[ht]
\centering
\scalebox{.8}{
\begin{minipage}{\linewidth}
\begin{subfigure}[t]{.49\linewidth}
\includegraphics[scale=.5]{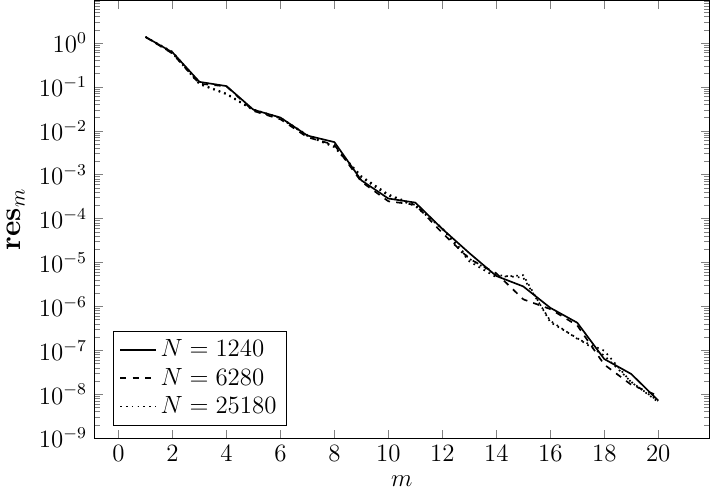}
\caption{\texttt{canyon\_particle}: SV-AAA}
\end{subfigure}\hfill
\begin{subfigure}[t]{.49\linewidth}
\includegraphics[scale=.5]{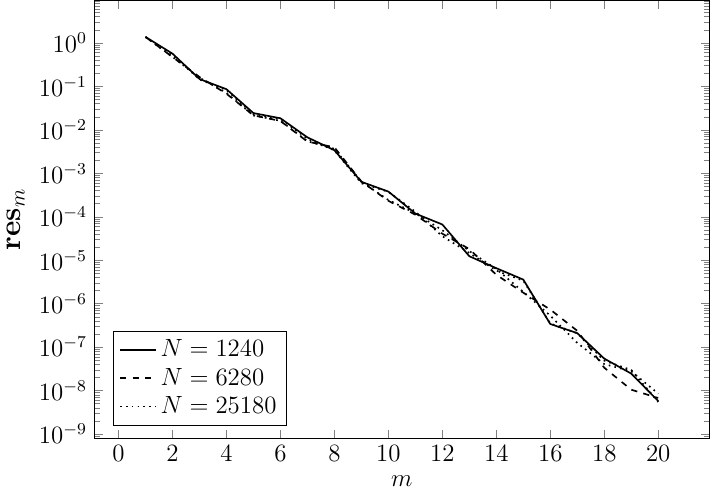}
\caption{\texttt{canyon\_particle}: QR-AAA}
\end{subfigure}

\vspace{\baselineskip}
\begin{subfigure}[t]{.49\linewidth}
\includegraphics[scale=.5]{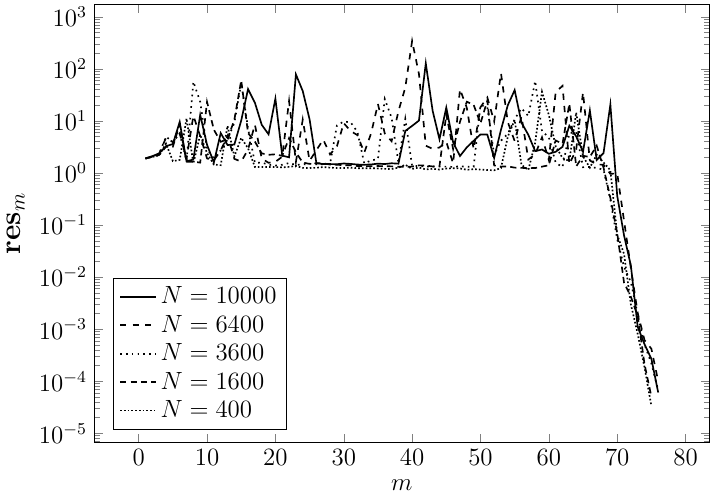}
\caption{\texttt{time\_delay3}: SV-AAA}
\end{subfigure}
\hfill
\begin{subfigure}[t]{.49\linewidth}
\includegraphics[scale=.5]{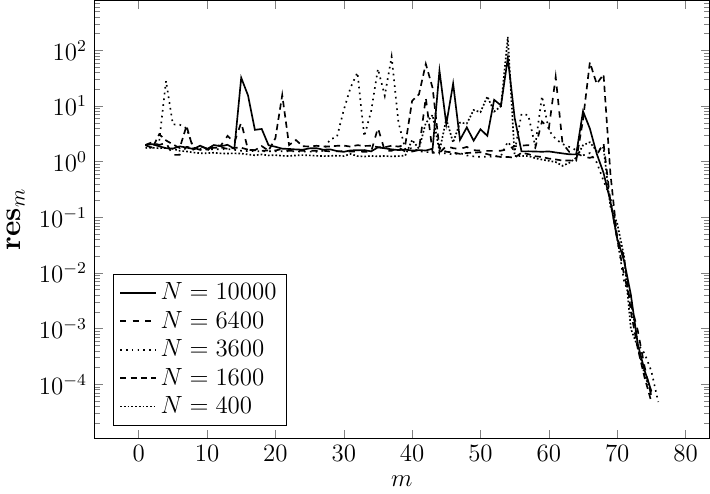}
\caption{\texttt{time\_delay3}: QR-AAA}
\end{subfigure}
\hfill
\end{minipage}
}
\caption{Residue $\textbf{res}_m$ for SV-AAA (left) and QR-AAA (right) over the degree $m$, for the second two selected problems.}\label{fig:errQRAAAPart2}
\end{figure}

We ran our experiments on an 8 core (16 virtual) AMD Ryzen 7 PRO 5850U x86\_64 CPU and report the timings and convergence behavior of both QR-AAA and SV-AAA. Each of the four chosen NLEVP problems comes in various sizes $N$, corresponding to the number of functions to be approximated.\\
The time (in seconds) needed for the QR approximation of $F\in\C^{Z\times N}$ is denoted by \texttt{t\_QR}, while the time for the set-valued AAA approximations for the $Q$ factor and the $F$ matrix are denoted by \texttt{t\_AAA\_Q} and \texttt{t\_AAA\_F} respectively. The timings were averaged over 10 experiments each, except for \texttt{t\_AAA\_F} for the largest instance (at $N=98442$) for problem $4$. The code for the four NLEVP QR-AAA experiments can be found in \cite{PQRAAAGit}.\\
We report the values of the absolute timings for \texttt{t\_QR}, \texttt{t\_AAA\_Q}, \texttt{t\_AAA\_F} in the tables in Figure~\ref{fig:tablesQRAAA}. When compared to the timings for SV-AAA, the timings for QR-AAA are almost negligible. Only at very large scales, basically when a single column of $Q$ no longer fits in cache, do we even register any cost for the QR decomposition at all. The only exception to this is problem $4$. We clearly see that \texttt{t\_AAA\_Q} does not increase with increasing $N$, in fact, it often \emph{decreases}. From the convergence plots in Figures~\ref{fig:errQRAAAPart1} and~\ref{fig:errQRAAAPart2}, we see why this is; the timings for the set-valued part of QR-AAA depend on the degree of approximation $m$, and the rank of $Q$, but not on $N$. The time \texttt{t\_QR} however does increase linearly with $N$, as expected, but since the QR factorization is fundamentally less computationally expensive than the SVD of the large Loewner matrices arising in the direct application of SV-AAA on $F$, we have that $\texttt{t\_QR}\ll\texttt{t\_AAA\_F}$.\\
Finally, we can also see from Figures~\ref{fig:errQRAAAPart1} and~\ref{fig:errQRAAAPart2} that the SV-AAA approximation and the QR-AAA approximation are essentially identical. This shows that QR-AAA is essentially equivalent to SV-AAA, only much faster.

\begin{figure}[ht]
\centering
\begin{subfigure}[t]{.47\linewidth}
\renewcommand{\arraystretch}{1.2}
\begin{tabular}{|c|c|c|c|}
\hline
$N$ & \texttt{t\_QR}&\texttt{t\_AAA\_Q}&\texttt{t\_AAA\_F}\\
\hline
$1240$ & $0.025$&$4.02\mathrm{e}{-3}$&$1.93$\\
\hline
$6280$ & $0.124$&$2.24\mathrm{e}{-3}$&$6.04$\\
\hline
$25180$ & $0.470$&$1.63\mathrm{e}{-3}$&$20.41$\\
\hline
\end{tabular}
\caption{\texttt{sandwich\_beam}}
\end{subfigure}
\hfill
\begin{subfigure}[t]{.47\linewidth}
\renewcommand{\arraystretch}{1.2}
\begin{tabular}{|c|c|c|c|}
\hline
$N$ & \texttt{t\_QR}&\texttt{t\_AAA\_Q}&\texttt{t\_AAA\_F}\\
\hline
$5546$ & $0.108$&$5.56\mathrm{e}{-3}$&$2.61122$\\
\hline
$29034$ & $0.527$&$1.26\mathrm{e}{-3}$&$13.7263$\\
\hline
$98442$ & $1.756$&$1.14\mathrm{e}{-3}$&$54.847$\\
\hline
\end{tabular}
\caption{\texttt{photonic\_crystal}}
\end{subfigure}\\
\vspace{\baselineskip}
\begin{subfigure}[b]{.47\linewidth}
\renewcommand{\arraystretch}{1.2}
\begin{tabular}{|c|c|c|c|}
\hline
$N$ & \texttt{t\_QR}&\texttt{t\_AAA\_Q}&\texttt{t\_AAA\_F}\\
\hline
$5971$ & $0.278$&$8.351\mathrm{e}{-2}$&$70.4021$\\
\hline
$7432$ & $0.381$&$8.39\mathrm{e}{-2}$&$86.9885$\\
\hline
$10157$ & $0.447$&$8.32\mathrm{e}{-2}$&$119.217$\\
\hline
$15121$ & $0.683$&$8.38\mathrm{e}{-2}$&$177.209$\\
\hline
\end{tabular}
\caption{\texttt{canyon\_particle}}
\end{subfigure}
\hfill
\begin{subfigure}[b]{.47\linewidth}
\renewcommand{\arraystretch}{1.2}
\begin{tabular}{|c|c|c|c|}
\hline
$N$ & \texttt{t\_QR}&\texttt{t\_AAA\_Q}&\texttt{t\_AAA\_F}\\
\hline
$400$ & $0.010$&$9.21$&$232.38$\\
\hline
$1600$ & $0.091$&$9.41$&$1033.47$\\
\hline
$3600$ & $0.241$&$9.51$&$2285.06$\\
\hline
$6400$ & $0.440$&$9.67$&$4431.65$\\
\hline
$10000$ & $0.701$&$9.41$&$6370.27$\\
\hline
\end{tabular}
\caption{\texttt{time\_delay3}}
\end{subfigure}
\caption{Tables of absolute timings for SV-AAA and QR-AAA, for the four NLEVP problems. All timings are reported in seconds.}\label{fig:tablesQRAAA}
\end{figure}

\begin{figure}[H]
\centering
\scalebox{.9}{
\begin{minipage}{\linewidth}
\begin{subfigure}[b]{.47\linewidth}
\includegraphics[scale=.5]{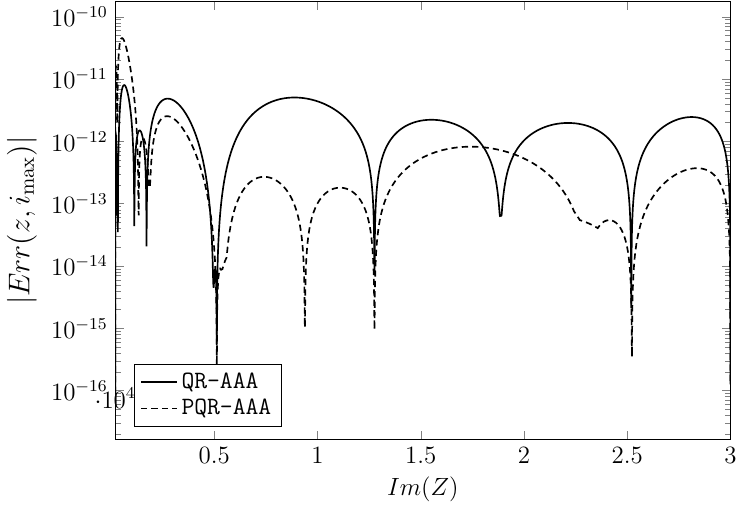}
\end{subfigure}
\hfill
\begin{subfigure}[b]{.47\linewidth}
\includegraphics[scale=.5]{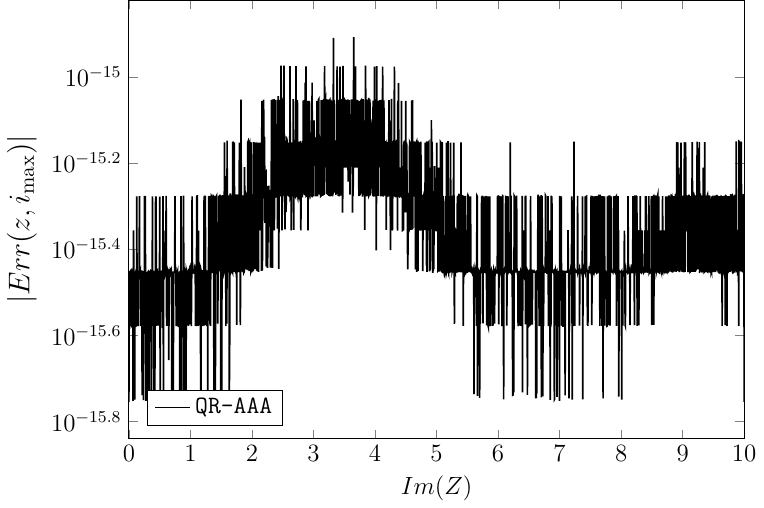}
\end{subfigure}
\\
\begin{subfigure}[b]{.47\linewidth}
\includegraphics[scale=.5]{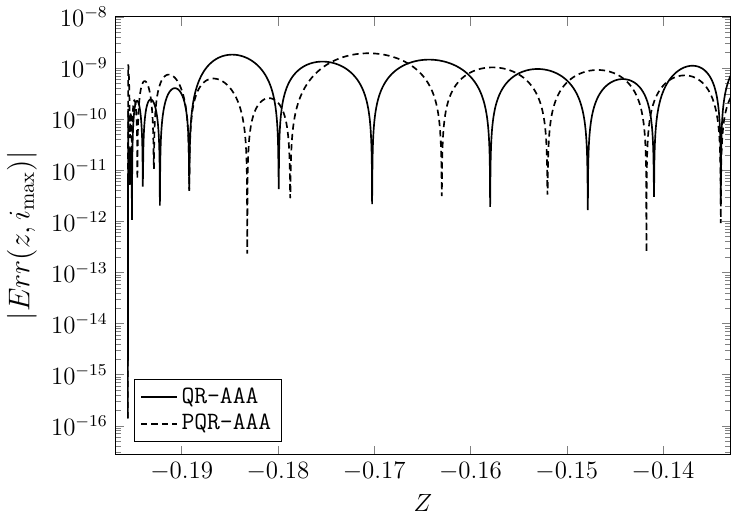}
\end{subfigure}
\hfill
\begin{subfigure}[b]{.47\linewidth}
\includegraphics[scale=.5]{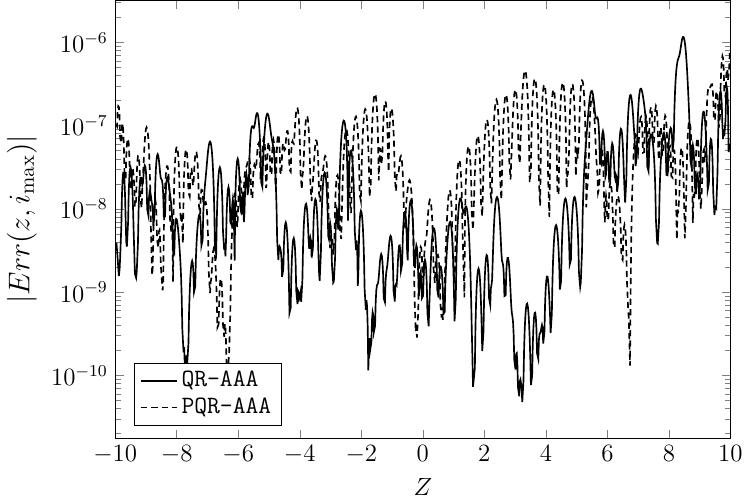}
\end{subfigure}
\end{minipage}
}
\caption{Relative $\|\cdot\|_{\infty}$-error over the domain of approximation for QR-AAA (full) and PQR-AAA (dashed) for the 4 NLEVP problems. From top to bottom and left to right: \texttt{sandwich\_beam}, \texttt{photonic\_crystal}, \texttt{canyon\_particle}, \texttt{time\_delay3}}\label{fig:errQRAAA}
\end{figure}

While in Figures~\ref{fig:errQRAAAPart1} and \ref{fig:errQRAAAPart2} we can see the convergence of the QR-AAA and SV-AAA residue (and hence also their final degree of approximation), in Figure~\ref{fig:errQRAAA} we plot the actual supremum norm error of our final QR-AAA and PQR-AAA approximations (using $4$ cores) over their respective domains. These errors were computed in relative $\|\cdot\|_{\infty}$-norm, sampled on a finer grid $|Z_{\text{test}}|=2513$. Note that, at $m=7$, the QR-AAA and PQR-AAA approximations are essentially exact for the \texttt{photonic\_crystal} problem. For QR-AAA, this was already apparent in Figure~\ref{fig:errQRAAAPart1}.
We see that, while PQR-AAA selects different support points, its accuracy is comparable. This supports the theoretical analysis from Section~\ref{sec:parallel}.
\subsection{Comparison to weighted AAA}
As discussed in Section~\ref{sec:intro}, for split form NLEVPs, our method is similar to weighted AAA (without refinement) from \cite{robustRatGuttel}. In this section we compare the degree and error obrained by the two methods, for a given tolerance. Note that weighted AAA exploits the split form structure, where QR-AAA does not, and as such weighted AAA is always significantly faster; its performance gain grows as the size of the matrices in the split form (correspondingly, $N$ in QR-AAA) grows. This may seem like a significant disadvantage to using QR-AAA, however, there are many situations (e.g., the application outlined in Section \ref{sec:BEM}) in which the problem is not in a (known) split form. We use the \texttt{MATLAB} implementation of weighted AAA available at \cite{MatlabRobust}. In the interest of brevity, we report our results only for the smallest problem sizes $N$. The results for larger problem sizes are similar. We compute for both approximations the relative $\|\cdot\|_{\infty}$-error over a more finely sampled grid $|Z_{test}|=2513$ corresponding to the same interval of approximation (e.g., $[200\imath,30000\imath]$ for the \texttt{sandwich\_beam} problem).
\begin{figure}[H]
\centering
\renewcommand{\arraystretch}{1.3}
\begin{tabular}{|c|c|c|c|c|}
\hline
Problem & $m_{QR}$&$m_{W}$&$\textbf{err}_{\infty,\text{QR}}$&$\textbf{err}_{\infty,\text{W}}$\\
\hline
\texttt{sandwich\_beam} ($N=1240$) & 7&6&$8.02\cdot10^{{-}12}$&$1.05\cdot10^{{-}10}$\\
\hline
\texttt{photonic\_crystal ($N=5546$)} & 7&5&$1.23\cdot10^{{-}15}$&$7.34\cdot10^{{-}10}$\\
\hline
\texttt{canyon\_particle} ($N=5971$) & 19&16&$3.5\cdot10^{{-}9}$&$7.21\cdot10^{{-}9}$\\
\hline
\texttt{time\_delay3} ($N=400$)& 84&85&$1.12\cdot 10^{{-}9}$&$3.86\cdot10^{{-}10}$\\
\hline
\end{tabular}
\caption{Comparison of QR-AAA and weighted AAA from \cite{robustRatGuttel} at given tolerance $10^{{-}8}$. Here $m_{\text{QR}}$ and $m_{\text{W}}$ denote number of support points in the approximation obtained by QR-AAA and weighted AAA respectively. Similarly, $\textbf{err}_{\infty,\text{QR}}$ and $\textbf{err}_{\infty,\text{W}}$ denote their respective relative $\|\cdot\|_{\infty}$-error over a finely sampled test grid in the domain of approximation.}
\label{fig:comparison1}
\end{figure}
We can see that, in terms of accuracy both methods perform similarly well. Both converge to well within the user-supplied tolerance. The $\|\cdot\|_{\infty}$-error is slightly better minimized over the entire interval by QR-AAA, at the cost of returning slightly more support points on average.
\subsection{BEM near-field compression using PQR-AAA}\label{sec:BEM}
\begin{figure}[ht]
\centering
\includegraphics[scale=.65]{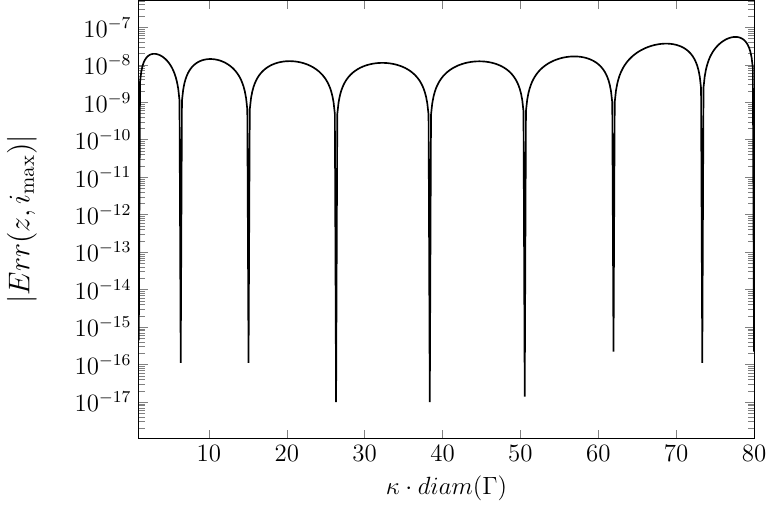}
\caption{Supremum norm error of the PQR-AAA approximation to $S(\kappa)$ over the near-field, varying with the dimensionless wave number $\kappa$}\label{fig:errplotNearField}
\end{figure}

In the Galerkin boundary element method for the Helmholtz equation (see, e.g., \cite{Sauter}, \cite{McLean}), the differential equation in a volume $\Omega$ is transformed into integral equations involving integral operators on the boundary $\partial\Omega$. The simplest such boundary integral operator, the single layer boundary operator, is defined by
$$\mathcal{S}(\kappa):H^{-1/2}(\partial\Omega)\to H^{1/2}(\partial\Omega):u\mapsto \mathcal{S}(\kappa)[u]$$
where
$$(\mathcal{S}(\kappa)[u])(\mathbf{x}) = \lim_{\mathbf{x}\to\partial\Omega}\int_{\partial\Omega}G(\mathbf{x},\mathbf{y};\kappa)u(\mathbf{y})dS_{\mathbf{y}}$$
with $G(\mathbf{x},\mathbf{y};\kappa)$ the Green's kernel at wave number $\kappa$. The operator $\mathcal{S}(\kappa)$ is discretized to $S(\kappa)$ using Galerkin discretization on a triangular surface mesh. This is a matrix function that is not in split form. In the nearfield, i.e., blocks of $S$ that correspond to clusters in $\partial\Omega$ that are not well-separated (see \cite{HMAT} for more details), costly specialized singular quadrature schemes must be used to ensure accuracy. We apply PQR-AAA to compress the near-field of $S(\kappa)$ as it varies over the wave number. This is an important task, as approximating the wavenumber dependence then allows us to avoid the costly near-field construction in the future. The matrix $F$ from equation~\ref{eq:defFmat} in this case has columns corresponding to $S(i,j,;Z)$, where $(i,j)$ is an index in the near-field DOFs, and $Z\subset[\kappa_{\min},\kappa_{\max}]$ is a discrete, sufficiently fine subset of the selected frequency range. For our experiment, we discretize $\mathcal{S}$ on a spherical grid obtained from a $6$-fold spherical refinement of an octahedral mesh using the software package \texttt{BEACHpack} (available at \cite{Beach}). We ran our experiment over 28 cores. The total number of functions to be approximated is $1000000$, which are distributed (roughly) equally among the cores. The required tolerance was set to $10^{\smin6}$. The (dimensionless\footnote{The dimensionless wavenumber is the actual wavenumber multiplied by the diameter of $\partial\Omega$}) wavenumber range was set to $[1.,80.]$, and discretized into $500$ equispaced points. The timings are split up into four constituents:
\begin{itemize}
\item \texttt{t\_F}: the maximal time needed to assemble the matrix $F_{\mu}$, over the processors $\mu\in\{1,\ldots,28\}$,
\item \texttt{t\_QR}: maximal time needed to compute the approximate pivoted Householder QR decomposition $F_{\mu}\approx Q_{\mu}R_{\mu}$ over the processors $\mu\in\{1,\ldots,28\}$,
\item \texttt{t\_AAA\_Q}: the maximal time needed to compute the SV-AAA approximation of $Q_{\mu}$ over the processors $\mu\in\{1,\ldots,28\}$.
\item \texttt{t\_fin}: the time needed to compute the final SV-AAA for $\uplus_{\mu}Q_{\mu}$
\end{itemize}
These timings, as well as the supremum norm error over $F$, are reported for our set-up in table~\ref{tab:timingsQRAAANearfield}. The error over the wavenumber of the final approximation is reported in Figure~\ref{fig:errplotNearField}. As can be seen from this figure, the final degree of rational approximation is $12$. The maximum error of the final approximant is well below the requested tolerance.
\begin{table}[ht]
\centering
\begin{tabular}{|c|c|c|c|c|}
\hline
\texttt{t\_F}&\texttt{t\_QR}&\texttt{t\_AAA\_Q}&\texttt{t\_fin}&$\textbf{Err}_{\infty}(F)$\\
\hline
$343.213$&$.926$&$.019$&$.259$&$6.64897\cdot10^{-8}$\\
\hline
\end{tabular}
\caption{Error and components of the PQR-AAA timings (in seconds) for the the near-field of $S(\kappa)$ defined on the sphere $S^2$.}\label{tab:timingsQRAAANearfield}
\end{table}

\section*{Acknowledgments}
The authors would like to extend their gratitude to Kobe Bruyninckx for his valuable help in writing the code for \texttt{BeachPack}, from which \texttt{pqr-aaa} heavily borrows. We would also like to thank the reviewers for their valuable feedback, which has greatly improved the presentation of this manuscript. This research was supported in part by KU Leuven IF project
C14/15/055 and FWO Research Foundation Flanders G0B7818N, as well as the National Science Foundation (DMS-2313434), and the Department of Energy ASCR (DE-SC0022251).

\bibliographystyle{acm}
\bibliography{bibQR}

\begin{thebibliography}{10}

\bibitem{AdcockOptimSampling}
{\sc Adcock, B., Platte, R.~B., and Shadrin, A.}
\newblock Optimal sampling rates for approximating analytic functions from pointwise samples.
\newblock {\em IMA J. Numer. Anal. 39}, 3 (2018), 1360--1390.

\bibitem{scalarRationalAAAAntoulas}
{\sc Antoulas, A.~C., and Anderson, B. D.~Q.}
\newblock On the scalar rational interpolation problem.
\newblock {\em {IMA} J. Math. Contr. Inf. 3}, 2 (1986), 61--88.

\bibitem{NLEVP1}
{\sc Betcke, T., Higham, N.~J., Mehrmann, V., Schr\"{o}der, C., and Tisseur, F.}
\newblock {NLEVP}: A collection of nonlinear eigenvalue problems.
\newblock {\em {ACM} Trans. Math. Softw. 39}, 2 (2013), 1--28.

\bibitem{NLEVPGit}
{\sc Betcke, T., Higham, N.~J., Mehrmann, V., Schr\"{o}der, C., and Tisseur, F.}
\newblock {NLEVP}: A collection of nonlinear eigenvalue problems.
\newblock https://www.github.com/ftisseur/nlevp, 2013.

\bibitem{BoydMockCheb}
{\sc Boyd, J.~P., and Xu, F.}
\newblock Divergence ({R}unge phenomenon) for least-squares polynomial approximation on an equispaced grid and mock {C}hebyshev subset interpolation.
\newblock {\em Appl. Math. Comput. 210}, 1 (2009), 158--168.

\bibitem{Beach}
{\sc Dirckx, S., and Bruyninckx, K.}
\newblock {BEACHpack}.
\newblock https://gitlab.kuleuven.be/numa/software/beachpack, 2023.

\bibitem{PQRAAAGit}
{\sc Dirckx, S., and Bruyninckx, K.}
\newblock pqraaa.
\newblock https://github.com/SimonDirckx/pqraaa, 2024.

\bibitem{Dirckx}
{\sc Dirckx, S., Huybrechs, D., and Meerbergen, K.}
\newblock Frequency extraction for {BEM} matrices arising from the 3{D} scalar {H}elmholtz equation.
\newblock {\em SIAM J. Sci. Comput. 44}, 5 (2022), B1282--B1311.

\bibitem{DongID}
{\sc Dong, Y., and Martinsson, P.-G.}
\newblock Simpler is better: a comparative study of randomized pivoting algorithms for cur and interpolative decompositions.

\bibitem{photonic}
{\sc Effenberger, C., Kressner, D., and Engström, C.}
\newblock Linearization techniques for band structure calculations in absorbing photonic crystals.
\newblock {\em International Journal for Numerical Methods in Engineering 89\/} (01 2012), 180 -- 191.

\bibitem{GuideSaadRatApprox}
{\sc El-Guide, M., Mi\c{e}dlar, A., and Saad, Y.}
\newblock A rational approximation method for solving acoustic nonlinear eigenvalue problems.
\newblock {\em Eng Anal Bound Elem 111\/} (2020), 44--54.

\bibitem{RatMinMax}
{\sc Filip, S.-I., Nakatsukasa, Y., Trefethen, L.~N., and Beckermann, B.}
\newblock Rational minimax approximation via adaptive barycentric representations.
\newblock {\em SIAM Journal on Scientific Computing 40}, 4 (2018), A2427--A2455.

\bibitem{RBFQR}
{\sc Fornberg, B., Larsson, E., and Flyer, N.}
\newblock Stable computations with {G}aussian radial basis functions.
\newblock {\em SIAM Journal on Scientific Computing 33}, 2 (2011), 869--892.

\bibitem{GolubAndVanLoan}
{\sc Golub, G.~H., and van Loan, C.~F.}
\newblock {\em Matrix Computations}, {F}ourth~ed.
\newblock JHU Press, 2013.

\bibitem{robustRatGuttel}
{\sc G\"{u}ttel, S., Negri~Porzio, G.~M., and Tisseur, F.}
\newblock Robust rational approximations of nonlinear eigenvalue problems.
\newblock {\em SIAM J. Sci. Comput. 44}, 4 (2022), A2439--A2463.

\bibitem{MatlabRobust}
{\sc Güttel, S., Negri~Porzio, G.~M., and Tisseur, F.}
\newblock Nep2rat.
\newblock https://github.com/Gmnp/nep2rat, 2020.

\bibitem{HMAT}
{\sc Hackbusch, W.}
\newblock {\em Hierarchical Matrices: Algorithms and Analysis}.
\newblock Springer Series in Computational Mathematics. Springer-Verlag, Berlin-Heidelberg, 2010.

\bibitem{NLEVP2}
{\sc Higham, N.~J., Negri~Porzio, G.~M., and Tisseur, F.}
\newblock An updated set of nonlinear eigenvalue problems.
\newblock Available at {https://api.semanticscholar.org/CorpusID:164833102}, 2019.

\bibitem{HochmanFastAAA}
{\sc Hochman, A.}
\newblock Fast{AAA}: A fast rational-function fitter.
\newblock In {\em 2017 IEEE 26th Conference on Electrical Performance of Electronic Packaging and Systems (EPEPS)\/} (2017), pp.~1--3.

\bibitem{AAAequi}
{\sc Huybrechs, D., and Trefethen, L.~N.}
\newblock {AAA} interpolation of equispaced data.
\newblock {\em BIT 63\/} (2023), 1--19.

\bibitem{KarlSVAAA}
{\sc Lietaert, P., Meerbergen, K., P{\'e}rez, J., and Vandereycken, B.}
\newblock {Automatic rational approximation and linearization of nonlinear eigenvalue problems}.
\newblock {\em IMA J. Numer. Anal. 42}, 2 (2022), 1087--1115.

\bibitem{McLean}
{\sc McLean, W.}
\newblock {\em Strongly Elliptic Systems and Boundary Integral Equations}.
\newblock Cambridge University Press, 2000.

\bibitem{AAA}
{\sc Nakatsukasa, Y., S\`{e}te, O., and Trefethen, L.~N.}
\newblock The {AAA} algorithm for rational approximation.
\newblock {\em SIAM J. Sci. Comput. 40}, 3 (2018), A1494--A1522.

\bibitem{Sauter}
{\sc Sauter, S., and Schwab, C.}
\newblock {\em Boundary Element Methods}, vol.~39 of {\em Springer Series in Computational Mathematics}.
\newblock Springer, Berlin, 2011.

\bibitem{SuBaiSVRat}
{\sc Su, Y., and Bai, Z.}
\newblock Solving rational eigenvalue problems via linearization.
\newblock {\em SIAM J. Matrix Anal. Appl. 32}, 1 (2011), 201--216.

\bibitem{MeerbergenBeam}
{\sc Van~Beeumen, R., Meerbergen, K., and Michiels, W.}
\newblock A rational {K}rylov method based on {H}ermite interpolation for nonlinear eigenvalue problems.
\newblock {\em SIAM J. Sci. Comput. 35}, 1 (2013), A327--A350.

\bibitem{CORK}
{\sc Van~Beeumen, R., Meerbergen, K., and Michiels, W.}
\newblock Compact rational {Krylov} methods for nonlinear eigenvalue problems.
\newblock {\em SIAM J. Matrix Anal. Appl. 36\/} (2015), 820--838.

\bibitem{MeerbergenCanyon}
{\sc Vandenberghe, W., Fischetti, M., Van~Beeumen, R., Meerbergen, K., Michiels, W., and Effenberger, C.}
\newblock Determining bound states in a semiconductor device with contacts using a nonlinear eigenvalue solver.
\newblock {\em J. Comput. Electron. 13\/} (07 2014), 753--762.

\bibitem{VoroninMartinsson}
{\sc Voronin, S., and Martinsson, P.-G.}
\newblock Efficient algorithms for cur and interpolative matrix decompositions.
\newblock {\em Advances in Computational Mathematics 43\/} (06 2017), 1--22.

\bibitem{RBFflat}
{\sc Wright, G.~B., and Fornberg, B.}
\newblock Stable computations with flat radial basis functions using vector-valued rational approximations.
\newblock {\em JCP 331\/} (2017), 137--156.

\end{thebibliography}

\cleardoublepage

\end{document}